\documentclass[11pt]{amsart}

\usepackage{amssymb}

\usepackage{enumerate}

\usepackage{graphicx}

\makeatletter
\@namedef{subjclassname@2020}{%
  \textup{2020} Mathematics Subject Classification}
\makeatother

\newtheorem{thm}{Theorem}[section]
\newtheorem{cor}[thm]{Corollary}
\newtheorem{lem}[thm]{Lemma}
\newtheorem{prop}[thm]{Proposition}

\theoremstyle{definition}

\newtheorem{rem}[thm]{Remark}

\newtheorem{fact}[thm]{Fact}

\numberwithin{equation}{section}

\begin{document}

\baselineskip=17pt

\title[Three results in linear dynamics]{Three results in linear dynamics}

\author[Mohammad Ansari]{Mohammad Ansari}

\date{}
\begin{abstract} In this article, first we show that the Fr$\acute{\textnormal{e}}$chet space $H(\Bbb D)$ cannot support strongly supercyclic weighted composition operators.
Then we compute the constant $\epsilon$ for weighted backward shifts on $\ell^p$ ($1\le p<\infty$) and $c_0$. This constant is used to find strongly hypercyclic scalar
multiples of non-invertible strongly supercyclic Banach space operators.
 Finally, we give an affirmative answer to a recent open question concerning supercyclic vectors.
\end{abstract}

\maketitle
\section{Introduction}\label{sec1}
Suppose $\mathcal X$ is an infinite-dimensional second countable Baire topological vector space (over $\Bbb C$) and $\mathcal S$ is a multiplicative semigroup of maps on $\mathcal X$ 
 (the binary operation is the composition of two maps). For a vector $x\in \mathcal X$ we put  $$\mathcal Sx=\{f(x): f\in \mathcal S\}.$$ A semigroup $\mathcal S$ is said to be
 {\it hypercyclic} if there is a nonzero vector $x\in \mathcal X$ such that $\overline{\mathcal Sx}=X$. In this case, $x$ is called a {\it hypercyclic vector} for $\mathcal S$. 
 The set of all hypercyclic vectors for $\mathcal S$ is denoted by $HC(\mathcal S)$. If $HC(\mathcal S)=\mathcal X\backslash\{0\}$, then $\mathcal S$ is called 
{\it hypertransitive}. 
\par A semigroup $\mathcal S$ is called {\it supercyclic} if the semigroup $$\Bbb C\mathcal S=\{bf: b\in \Bbb C, f\in \mathcal S\}$$ is hypercyclic. 
A vector $x\in \mathcal X$ is called a {\it supercyclic vector} for $\mathcal S$ whenever it is a hypercyclic vector for $\Bbb C\mathcal S$. 
The set of all supercyclic vectors for $\mathcal S$ is denoted by $SC(\mathcal S)$.
\par We say that $\mathcal S$ is {\it topologically transitive} if, for any pair of nonempty open sets $U, V\subseteq \mathcal X$,
there is some $f\in \mathcal S$ such that $$f(U)\cap V\neq \emptyset.$$ Finally, $\mathcal S$ is said to be {\it strongly topologically transitive} if, for any nonempty open
subset $U$ of $\mathcal X$, we have that $$\mathcal X\backslash\{0\}\subseteq \bigcup_{f\in \mathcal S}f(U).$$ It is clear that strong topological transitivity implies topological transitivity.
\par One can easily verify that a semigroup $\mathcal S$ is hypertransitive if and only if $$\mathcal X\backslash\{0\}\subseteq \bigcup_{f\in \mathcal S}f^{-1}(U)$$ for every nonempty open subset 
$U$ of $\mathcal X$. Now assume that all members of $\mathcal S$ are bijections and put $$S^{-1}=\{f^{-1}: f\in \mathcal S\}.$$ Then the following fact is obvious:
\begin{fact} Let $\mathcal S$ be a semigroup in which every map is a bijection. Then $\mathcal S$ is strongly topologically transitive
 if and only if $\mathcal S^{-1}$ is hypertransitive.
\end{fact}
\par Recall that by $L(\mathcal X)$ we mean the set of all continuous linear operators on $\mathcal X$. From now on, by an operator, we mean a continuous linear operator.
The most particular case among all multiplicative operator semigroups is the semigroup generated by a single operator.
A continuous linear operator $T$ on $\mathcal X$ is said to be hypercyclic (resp. supercyclic, topologically transitive)  
if the semigroup $$\mathcal S=\{T^n: n\in \Bbb N_0\}$$ is hypercyclic (resp. supercyclic, topologically transitive). 
Here $\Bbb N_0$ is the set of all non-negative integers and $T^0=I$, the identity operator on $\mathcal X$.
We also say that $x$ is a hypercyclic (resp. supercyclic) vector for $T$ when it is a hypercyclic (resp. supercyclic) vector for $\mathcal S$, and we write $HC(T)$ (resp. $SC(T)$)
for $HC(\mathcal S)$ (resp. $SC(\mathcal S)$). Notice that (since $\mathcal X$ has no isolated points) if $x$ is a hypercyclic (resp. supercyclic) vector for $T$,
then so is every $T^px$ ($p\ge 1$) (resp. $bT^px$ ($0\neq b\in \Bbb C$, $p\ge 1$)), and this shows that $HC(T)$ (resp. $SC(T)$) is dense in $\mathcal X$.
\par Two well-known books which are always advised to study hypercyclic and supercyclic operators are \cite{bm,gp}.
\par We say that $T$ is {\it strongly hypercyclic} (resp. {\it strongly supercyclic}) if the semigroup $\mathcal S=\{T^n: n\in \Bbb N_0\}$ 
(resp. $\mathcal S=\{bT^n: b\in \Bbb C, n\in \Bbb N_0\}$) is
strongly topologically transitive. In other words, $T$ is strongly hypercyclic (resp. strongly supercyclic)
 if $$\mathcal X\backslash\{0\}\subseteq \bigcup_{n\in \Bbb N_0}T^n(U)\;\; \textnormal{\big{(}resp.}\;\; \mathcal X\backslash\{0\}\subseteq \bigcup_{b\in \Bbb C, n\in \Bbb N_0}bT^n(U)\big{)}$$
 for any nonempty open subset $U$ of $\mathcal X$. 
 \par It is clear that strong hypercyclicity (resp. strong supercyclicity) is stronger than hypercyclicity (resp. supercyclicity). 
 In fact, this claim is verified by taking a look at \cite[Theorem $1.2$]{bm} and its following remark (resp. \cite[Theorem $1.12$]{bm}).
 The reader is referred to \cite{a, ahkh} for more information on strong topological transitivity, strong hypercyclicity, and strong supercyclicity. 
 \par In Section $2$, we show that there is no strongly supercyclic weighted composition operator on the function space $H(\Bbb D)$. 
 In Section $3$, we compute the constant $\epsilon$
 for weighted backward shift operators on $\ell^p$ ($1\le p<\infty$) and $c_0$. In the last section, we give an affirmative answer to a recent open question about supercyclic vectors.
\section{Weighted composition operators on $H(\Bbb D)$}\label{sec2}
 It is proved in \cite{ahkh} that no automorphism invariant weighted Hardy
space $H^2(\beta)$ can support strongly supercyclic weighted composition operators. Some particular cases of theses spaces
are the classical Hardy space $H^2(\Bbb D)$ and the Bergman and Dirichlet spaces.
\par In the following theorem, we show that the same assertion is true
for the Fr$\acute{\textnormal{e}}$chet space $H(\Bbb D)$.
Recall that the space $H(\Bbb D)$ is the set of all analytic functions on $\Bbb D$, equipped with the compact-open
topology. It is worth mentioning that $H(\Bbb D)$ supports hypercyclic (and hence supercyclic) weighted composition operators \cite{yr}.
\par Assume that $w, \phi \in H(\Bbb D)$ with $\phi (\Bbb D)\subseteq \Bbb D$. The weighted composition operator
$C_{w,\phi}$ on $H(\Bbb D)$ is defined by $$C_{w,\phi}(f) = w(f\circ \phi),\;\textnormal{i.e.,}\;C_{w,\phi}(f)(z) = w(z)f(\phi(z)).$$
Meanwhile, it is easy to see that the $n$th ($n\ge 2$) iterate of $C_{w,\phi}$ is defined by
$$C^n_{w,\phi}(f) = w(w\circ \phi)\cdots (w\circ \phi_{n-1})(f\circ \phi_n)$$
for every $f\in H(\Bbb D)$, where $\phi_k$ is the $k$th iterate of $\phi$.
\begin{thm} The Fr$\acute{\textnormal{e}}$chet space $H(\Bbb D)$ cannot support strongly supercyclic weighted composition operators.
\end{thm}
\begin{proof} To get a contradiction, suppose that $C_{w,\phi}$ is a strongly supercyclic operator on $H(\Bbb D)$. 
Then we claim that $w(z)\neq 0$ for all $z\in \Bbb D$. To prove this claim, assume that $g$ is
the constant function $g(z)=1$ ($z\in \Bbb D$). Then, by the definition of strong supercyclicity, we must have
\begin{align}
g\in\bigcup_{b\in \Bbb C, n\in \Bbb N_0}bC^n_{w,\phi}(U)
\end{align}
 for any nonempty open subset $U$ of $H(\Bbb D)$.
Hence, there are some $f\in U$, $b\in \Bbb C$, and $n\in \Bbb N_0$ such that $g=bC^n_{w,\phi}f$. Without loss of generality, we can assume
that $n\ge 1$. Indeed, if $n=0$ then we can replace $U$ by the open set $V=U\backslash \Bbb Cg$ in $(2.1)$. This shows that, for all $z\in \Bbb D$, we have that
$$1=bw(z)w(\phi (z))\cdots w(\phi_{n-1}(z))f(\phi_{n}
(z)),$$ 
which proves that $w(z)\neq 0$ for all $z\in \Bbb D$. 
\par On the other hand, it is easy to see that $\phi$
cannot be a constant map. Indeed, if $\phi$ is constant, then we
have $$\{bC^n_{w,\phi}f : b\in \Bbb C, n\in \Bbb N_0\} \subseteq \Bbb C{w}$$
 for any $f\in H(\Bbb D)$. This contradicts our assumption that $C_{w,\phi}$ is strongly supercyclic. 
\par Hence $C_{w,\phi}$ is injective and so it is
a bijection by \cite[Proposition $3.5$]{ahkh}. Then, by \cite[Theorem 2.2]{b}, $\phi$ is an automorphism of $\Bbb D$ (the proof of that theorem only uses the assumption that $C_{w,\phi}$ is
a bijection). Now it is easily seen that $C_{w,\phi}C_{(1/w)\circ \phi^{-1},\phi^{-1}} = I$, and so we must have
$C^{-1}_{w,\phi} = C_{(1/w)\circ \phi^{-1},\phi^{-1}}$.
Note that $C^{-1}_{w,\phi}$ is a linear map on $H(\Bbb D)$ which may not necessarily
be continuous. Let us put $\psi = (1/w)\circ \phi^{-1}$ and $\rho=\phi^{-1}$
to simply write $C^{-1}_{w,\phi} = C_{\psi,\rho}$. 
\par Now, in view of the definition of strong supercyclicity, the semigroup
$$\mathcal S = \{bC^n_{w,\phi} : b\in \Bbb C\backslash\{0\}, n\in \Bbb N_0\}$$
is strongly topologically transitive (if $\mathcal S$ is strongly topologically transitive, then so is $\mathcal S\backslash\{0\}$), and hence, the semigroup
$$\mathcal S^{-1} = \{bC^n_{\psi,\rho} : b\in \Bbb C\backslash\{0\}, n\in \Bbb N_0\}$$
is hypertransitive by Fact $1.1$. So, every $f\in H(\Bbb D)$ which
is not identically zero is a hypercyclic vector for $\mathcal S^{-1}$. 
\par We consider the two
possible cases for the automorphism $\rho$:\\
$(1)$ There is some $a\in \Bbb D$ such that $\rho(a) = a$. Then the function $f(z) = z-a$ ($z\in \Bbb D$) is a hypercyclic
vector for $\mathcal S^{-1}$, and hence, there are sequences $(b_k)_k$ in $\Bbb C$ and $(n_k)_k$ in $\Bbb N_0$ such that
$b_kC^{n_k}_{\psi,\rho}f \to g$ where $g(z) = 1$ ($z\in \Bbb D$). Thus, at $z = a$ we must have
$$0 = b_k\psi(a)\psi(\rho(a))\cdots \psi(\rho_{n_k-1}(a))f(\rho_{n_k}(a))\to 1$$
which is not true.\\
$(2)$ The map $\rho$ has no fixed point in $\Bbb D$. Then, by the well-known Denjoy-Wolff theorem, there
is some $a\in \Bbb T$ (the unit circle $\partial \Bbb D$) such that $\rho_n\to a$ uniformly on compact subsets of $\Bbb D$. If we put
$A = \{\rho_n(0) : n\in \Bbb N\}$, then there is some $f\in H(\Bbb D)$ such that $Z(f) = A$ \cite[Theorem $15.11$]{r},
 where $Z(f)$ is the set of all zeros of $f$. Now $f$ is a hypercyclic
vector for $\mathcal S^{-1}$ and so, for the constant function $g(z) = 1$ ($z\in \Bbb D$), there exist sequences $(b_k)_k$ in $\Bbb C$ and $(n_k)_k$ in $\Bbb N_0$ such that 
$b_kC^{n_k}_{\psi,\rho}f \to g$. Then we have that
$$0 = b_k\psi(0)\psi(\rho(0))\cdots \psi(\rho_{n_k-1}(0))f(\rho_{n_k}(0))\to 1,$$ which is impossible.
\par Therefore, the assumption of strong supercyclicity of $C_{w,\phi}$ cannot be true and
the proof is complete.
\end{proof}
\section{The constant $\epsilon$ for weighted backward shifts}\label{sec3}
Let $\mathcal X$ be a Banach space and $T$ be a bounded linear operator on $\mathcal X$. In \cite{ahkh}, the constant $\epsilon(T)$ is defined by
$$\epsilon(T)=\inf\{\|y\|: y\in \mathcal X\backslash T(B)\},$$ where $B$ is the open unit ball of $\mathcal X$. Then it is proved that: 
\begin{thm} \cite[Theorem $4.2$]{ahkh} Assume that $T\in B(\mathcal X)$ is not invertible. Then $T$ is surjective with dense generalized
 kernel if and only if $cT$ is strongly
 hypercyclic for all $c\in \Bbb C$ with $|c|>1/\epsilon(T)$.
 \end{thm}
 \begin{cor} \cite[Corollary $4.3$]{ahkh} If $T$ is not invertible then
 $T$ is strongly supercyclic if and only if it is surjective and has dense generalized kernel. 
 \end{cor}
 Thus, the evaluation of $\epsilon(T)$ is
 important to find strongly hypercyclic scalar multiples of non-invertible strongly supercyclic Banach space operators. 
 \begin{rem} It is worth mentioning that, while in view of Theorem $3.1$ and Corollary $3.2$, non-invertible strongly supercyclic and
strongly hypercyclic Banach space operators are scalar multiples of one another, there are (non-invertible) supercyclic Banach space operators whose scalar multiples are never
hypercyclic \cite[Example $1.15$]{bm}. 
\end{rem}
 In the next result, we compute
 the constant $\epsilon$ for weighted backward shifts on $\ell^p$ ($1\le p<\infty$) and $c_0$.
Recall that $\ell^p$ ($1\le p<\infty$) is the space of all
complex sequences $(a_n)_{n\ge 0}$ satisfying
$\sum_{n=0}^{\infty}|a_n|^p<\infty$. The norm of $x=(a_n)_n$ in
$\ell^p$ is defined by $\|x\|=(\sum_{n=0}^{\infty}|a_n|^p)^{1/p}$.
The space $c_0$ is comprised of all sequences $(a_n)_{n\ge 0}$ in
$\Bbb C$ such that $a_n\to 0$, and the norm of $x=(a_n)_n$ in
$c_0$ is defined by $\|x\|=\sup_{n\ge 0}|a_n|$.
\par The weighted backward shift $B_W$ on $\mathcal X=\ell^p$ ($1\le p< \infty$)
or $c_0$ is defined by $B_W(e_0)=0$, and $B_W(e_n)=w_ne_{n-1}$
($n=1,2,3,\cdots$), where $W=(w_n)_{n\ge 1}$ is a bounded sequence of
positive numbers and $(e_n)_{n\ge 0}$ is the canonical basis of $\mathcal X$. 
One can readily verify that $B_W$ is
surjective if and only if $(w_n)_n$ is bounded-away from zero,
i.e., $\inf_{n\ge 1}w_n>0$. 
\begin{prop} Suppose that $B_W$ is a weighted backward shift on $\mathcal X=\ell^p$ ($1\le p< \infty$)
or $c_0$ with the weight sequence
$(w_n)_n$. Then we have that $\epsilon(B_W)=\inf_{n\ge 1}w_n$.
\end{prop}
\begin{proof} Put $\inf_{n\ge 1}w_n=r$. If $r=0$ then there exists a
strictly increasing sequence $(n_k)_k$ of positive integers such
that $w_{n_k}\to 0$ as $k\to \infty$. For each $k\ge 1$, put
$x_k=(0, 0, \cdots, 1, 0, 0, \cdots)$, where the number $1$ has
been set at the $n_k$-th position (remember that the position
numbering starts with zero). Now, if we put $y_k=B_Wx_k$ ($k\ge
1$), then it is easy to see that $y_k\in \mathcal X\backslash B_W(B)$. In
fact, if $z=(b_i)_{i\ge 0}$ and $y_k=B_Wz$ for a fixed $k\ge 1$,
then the equality $y_k=B_Wx_k$ shows that $b_{n_k}=1$, $b_i=0$ for
all $0< i\neq n_k$, and $b_0$ could be any complex number, and hence,
$\|z\|\ge 1$. Now, since $\|y_k\|=w_{n_k}\to 0$, we conclude that
$\epsilon(B_W)=0$.
 \par Now assume that $r>0$. Then $B_W$ is surjective. Let $y\in \mathcal X\backslash B_W(B)$ be an arbitrary
vector. Then there exists a vector $x=(a_0,a_1,\cdots)\notin B$
such that $y=B_Wx=(w_1a_1,w_2a_2,\cdots)$. If we set
$\hat{x}=(0,a_1,a_2,\cdots)$ then we must have $\|\hat{x}\|\ge 1$,
because otherwise, $y=B_Wx=B_W\hat{x}\in B_W(B)$ which contradicts
our assumption. Then $\|y\|\ge r\|\hat{x}\|\ge r$, and hence,
$\epsilon(B_W)\ge r$.
\par To complete our proof, we show that $\epsilon(B_W)\le
r$. Let $\delta>0$ be an arbitrary number. Then there is some
$k\ge 1$ such that $w_k<r+\delta$. Choose the vector
$x=(a_i)_{i\ge 0}$ for which $a_k=(r+\delta)w^{-1}_k$ and $a_i=0$
for all $0\le i\neq k$, and let $y=B_Wx$. Then $\|x\|=a_k>1$ and
meanwhile, we claim that $y\notin B_W(B)$. Indeed, if
$z=(b_i)_{i\ge 0}$ and $y=B_Wz$ then we have that $b_k=a_k$, $b_i=0$
for all $0< i\neq k$, and $b_0$ could be any complex number. Thus,
$\|z\|>1$. Now, since $y\notin B_W(B)$ and $\|y\|=r+\delta$, we
have $\epsilon(B_W)\le r+\delta$. Finally, since $\delta>0$ was
arbitrary, we conclude that $\epsilon(B_W)\le r$.
\end{proof}
\section{A positive answer to an open question}\label{sec4}
 We finish this paper by giving an affirmative answer to a recently asked open question concerning supercyclic vectors.
 It should be mentioned that this section has appeared as an arXiv preprint \cite{a1}.
\par For a subset $M$ of $\mathcal X$, by $\bigvee M$ we mean the closed linear span of $M$, i.e., $\bigvee M=\overline{\text{span}M}$.
 It is clear that if $x$ is a supercyclic vector 
for an operator $T$ on $\mathcal X$, then $$\bigvee \{T^nx: n\in \Bbb N_0\}=\mathcal X.$$ A natural question which may be asked here is that whether 
there is a strictly increasing sequence $(n_k)_k$ of positive integers such that $$\bigvee \{T^{n_k}x: k\ge 1\}\neq \mathcal X.$$
 \par In their recently published paper, Faghih-Ahmadi and Hedayatian \cite{fh} have proved the following interesting result. Recall that for a normed linear space $\mathcal X$, we often
 write $B(\mathcal X)$ instead of $L(\mathcal X)$.
\begin{thm}[Theorem $1$ of \cite{fh}] Let $\mathcal H$ be an infinite-dimensional Hilbert space. If $x$ is a supercyclic vector for $T\in B(\mathcal H)$, then
there is a (strictly increasing) sequence $(n_k)_k$ of positive integers such that $\bigvee \{T^{n_k}x : k\ge 1\}\neq \mathcal H$.
\end{thm}
Then they have asked whether the assertion is true for locally convex spaces or at least for Banach spaces \cite[Question $1$]{fh}. 
 \par We show that the assertion is true for normed linear spaces. To prove our result, we use Lemma $4.3$ which is an analogue of the following lemma.
\begin{lem}[Lemma $2.3$ of \cite{bkk}] Let $\mathcal A$ be a dense subset of an infinite-dimensional Banach space $\mathcal X$ and $e$ be a fixed element with $\|e\|>1$. 
Then, for every finite-dimensional subspace $Y\subset \mathcal X$ with $\text{dist}(e, Y )>1$, for every $\epsilon >0$ and $y\in Y$, there is an $a\in \mathcal A$ such that
$\|y-a\|<\epsilon$ and $\text{dist}(e, \text{span}\{Y, a\})>1$.
\end{lem}
 The proof of the above lemma shows that it can also be stated for infinite-dimensional normed linear spaces. 
 In fact, the lemma is proved by using
 Lemma $2.2$ of \cite{bkk} whose proof uses the fact that the weak closure of the unit sphere is the unit ball.
 Thus, we can give the following modified version.
\begin{lem} Let $\mathcal A$ be a dense subset of an infinite-dimensional normed linear space $\mathcal X$ and $e\in \mathcal X$ be a fixed element with $\|e\|>1$. 
Then, for every finite-dimensional subspace $Y$ of $\mathcal X$ with $\text{dist}(e,Y)>1$, there is an $a\in \mathcal A$ such that
$\text{dist}(e, \text{span}\{Y, a\})>1$.
\end{lem}
Now we are ready to answer Question $1$ of \cite{fh} for normed linear spaces.
\begin{thm} Let $\mathcal X$ be an infinite-dimensional normed linear space. If $x$ is a supercyclic vector for $T\in B(\mathcal X)$ then
there is a strictly increasing sequence $(n_k)_k$ of positive integers such that $\bigvee \{T^{n_k}x : k\ge 1\}\neq \mathcal X$.
\end{thm}
\begin{proof} It is clear that every nonzero scalar multiple of $x$ is also a supercyclic vector for $T$.
On the other hand, since $x$ and $Tx$ are linearly independent vectors, we have that $\text{dist}(x,\text{span}\{Tx\})>0$.
Thus, without loss of generality, we can assume that $\|x\|>1$ and $\text{dist}(x,\text{span}\{Tx\})>1$. Now, if
 we put $$\mathcal A=\{cT^nx: c\in \Bbb C, n>1\}, Y=\text{span}\{Tx\},$$ and $e=x$, then 
$\mathcal A$ is dense in $\mathcal X$ because $$\mathcal A=\{cT^nx: c\in \Bbb C, n\ge 0\}\backslash \{cT^nx: c\in \Bbb C, n=0, 1\}$$ and
$\{cT^nx: c\in \Bbb C, n=0, 1\}$ is nowhere dense in $\mathcal X$, and hence,
in view of Lemma $4.3$, there is
some $a=cT^{n_2}x\in \mathcal A$ such that $$\text{dist}(x,\text{span}\{Tx,T^{n_2}x\})=\text{dist}(x, \text{span}\{Y,a\})>1.$$ 
Now let $$\mathcal A_2=\{cT^nx: c\in \Bbb C, n> n_2\}, Y_2=\text{span}\{Tx,T^{n_2}x\}.$$ Then, by Lemma $4.3$,
there is some $a_2=c_2T^{n_3}x\in \mathcal A_2$ such that $$\text{dist}(x,\text{span}\{Tx,T^{n_2}x,T^{n_3}x\})=\text{dist}(x,\text{span}\{Y_2,a_2\})>1.$$ 
By continuing this construction, assume that for some $k\ge 2$, the dense set $\mathcal A_k$ and the finite-dimensional subspace $Y_k$ have been presented
 and (by using Lemma $4.3$) we have found an element $a_k=c_kT^{n_{k+1}}x\in \mathcal A_k$ such that
 $$\text{dist}(x,\text{span}\{Tx,T^{n_2}x,\cdots, T^{n_{k+1}}x\})=\text{dist}(x,\text{span}\{Y_k,a_k\})>1.$$ 
Then we put $$\mathcal A_{k+1}=\{cT^nx: c\in \Bbb C, n>n_{k+1}\}, Y_{k+1}=\text{span}\{Tx, T^{n_2}x, \cdots, T^{n_{k+1}}x\}.$$
Again, by Lemma $4.3$, there is some $a_{k+1}=c_{k+1}T^{n_{k+2}}x\in \mathcal A_{k+1}$ such that 
$$\text{dist}(x,\text{span}\{Tx,T^{n_2}x,\cdots, T^{n_{k+2}}x\})=\text{dist}(x,\text{span}\{Y_{k+1},a_{k+1}\})>1.$$ 
This inductive procedure gives a strictly increasing sequence $(n_k)_k$ (with $n_1=1$) and it is easily seen that $$\text{dist}(x,\bigvee\{T^{n_k}x: k\ge 1\})\ge 1.$$
Indeed, let $M=\text{span}\{T^{n_k}x: k\ge 1\}$ and suppose, to get a contradiction, that $\text{dist}(x,\overline{M})< 1$. Then there is some $y\in M$ such that
$\text{dist}(x,y)<1$. But it is clear that $y\in \text{span}\{Y_k, a_k\}$ for some $k\ge 1$ and we have already seen that $\text{dist} (x,\text{span}\{Y_k,a_k\})>1$.
Therefore, the assumption $\text{dist}(x,\overline{M})< 1$ cannot be true and this shows that $x\notin \bigvee\{T^{n_k}x: k\ge 1\}$.
\end{proof}
We need to mention that the authors in \cite{bb} have also answered the above-mentioned open question independently.
The interested readers are invited to investigate Question $1$ of \cite{fh} for operators on locally convex spaces.\\
{\bf Acknowledgments.} The author is very grateful to the reviewers for their hints and suggestions which improved the quality
of this paper.

\vspace*{.5cm}
\hspace*{.35cm} Department of Mathematics\\
\hspace*{.35cm} Azad University of Gachsaran\\
\hspace*{.35cm} Gachsaran, Iran\\
\hspace*{.35cm} {\it E-mail address}: \email{\small
ansari.moh@gmail.com}, \email{\small mo.ansari@iau.ac.ir}
\end{document}